\documentclass[12pt,reqno]{amsart}
\usepackage{amssymb,delarray}
\usepackage{amsfonts}
\usepackage{epsfig}
\usepackage[all]{xy}

\makeindex{}

\def\edvo{\rule {6pt}{6pt}}

\newtheorem{thm}{Theorem}
\newtheorem{lem}{Lemma}

\begin{document}
\baselineskip=14.pt plus 2pt 

\title[
]{On some quotient groups of hyperbolic groups}
\author[]{O.V.~Kulikova}

\address{\newline O.V.Kulikova \newline
Lomonosov Moscow State University, \newline Moscow Center for Fundamental and Applied Mathematics \newline (Moscow,  Russian Federation) }
\email{olga.kulikova@mail.ru}


\thanks{The work was supported by the Russian Science
Foundation, project no. 22-11-00075.}

\maketitle
\begin{abstract}
  This paper describes some generalizations of the results presented in the book "Geometry of defining Relations in Groups"\, of A.Yu.Ol'shanskii to the case of non-cyclic torsion-free hyperbolic groups. In particular, it is proved that for every non-cyclic torsion-free hyperbolic group, there exists a non-Abelian torsion-free quotient group in which all proper subgroups are cyclic, and the intersection of any two of them is not trivial.
\end{abstract}


  \setcounter{tocdepth}{2}
\def\st{{\sf st}}


\section*{ Introduction.}

Similar to free groups, hyperbolic groups have many different homomorphic images.

In the article \cite{olsh91}, A.Yu.Ol'shanskii noted that the constructions in \cite{olsh91} open up the possibility to generalize many of the results presented in his book \cite{olsh}
to the case of torsion-free hyperbolic groups. This paper describes some such generalizations of Theorems 31.1-31.4, 31.7, 31.8
from \cite{olsh}.
For example, in Theorem 2 of \cite{ashmanovOlsh} (= Theorem 31.4 of \cite{olsh}), I.S.Ashmanov and A.Yu.Ol'shanskii strengthened the example of S.I. Adyan \cite{adyan71, adyan75} and
proved the existence of a non-Abelian torsion-free group, in which all proper subgroups are cyclic, and the intersection of any two of them is not trivial.
In this paper, using the limit group constructed in \cite{semenov},  we prove that
for every non-cyclic torsion-free hyperbolic group, there exists a non-Abelian torsion-free quotient group, in which all proper subgroups are cyclic, and the intersection of any two of them is not trivial. Using this quotient, we get the generalization of Theorem 31.7 from \cite{olsh}:
for every prime $p$, for every non-cyclic torsion-free hyperbolic  group, there exists a non-Abelian quotient group $K$ of finite period, in which all proper subgroups are cyclic, and the Sylow $p$-subgroup $P$ is central, but $P$ is not a direct multiplier of $K$.

 In \cite{olsh93}, for every non-cyclic torsion-free hyperbolic group, A.Yu.Ol'shanskii constructed a non-Abelian quotient group, all of whose proper subgroups are cyclic, and the intersection of any two different maximal proper subgroups is trivial.
In this paper, using this quotient group, we obtain that for every non-cyclic torsion-free hyperbolic  group and for every $k\in\mathbb{N}$, there exists a non-Abelian torsion-free quotient group $\mathbf{A}$ in which all maximal proper subgroups are free Abelian groups of rank $k$ and the intersection of any two different maximal proper subgroups coincides with the center of $\mathbf{A}$, which is a free Abelian group of rank $k-1$, and no nontrivial subgroup of the center is a direct multiplier of $\mathbf{A}$.


 The proofs in this paper rely heavily on the proofs and results from \cite{olsh, ashmanovOlsh,semenov,olsh91,olsh93}.

The author thanks A.Yu. Ol'shanskii for the problem statement, the proposed method of solution and useful discussions during the writing of this article.

 \section{Designations}\label{S:df}

Let $G$ be a group with a generating set $\mathcal{A}$ and let $\mathcal{O}$ be the set of all relators of $G$.

 Let $\mathcal{R}=\bigcup_{i=1}^{\infty}\mathcal{R}_i$, where $\mathcal{R}_i$ is a symmetrized set of words in the alphabet $\mathcal{A}^{\pm 1}$, $\mathcal{R}_i\cap\mathcal{R}_j=\emptyset$ for $i\neq j$, some $\mathcal{R}_i$ may be empty, and $\mathcal{R}\cap\mathcal{O}=\emptyset$.

For every $k=1,2,\ldots$, define a group $G_k=\langle \mathcal{A}\|\mathcal{O}\cup \bigcup_{i=1}^{k}\mathcal{R}_i\rangle$.
 The presentation
 \begin{equation} \label{infPresent}
G(\infty) = \langle \mathcal{A}\|\mathcal{O}\cup \bigcup_{i=1}^{\infty}\mathcal{R}_i\rangle
\end{equation}
of a group $G(\infty)$ with such partition of the defining relation set will be called {\it  $G$-graded}.

By $\Phi$, we denote the canonical homomorphism from the free group $F=F(\mathcal{A})$ to the group $G=\langle \mathcal{A} \| \mathcal{O} \rangle$ with kernel $N$. Then $$G(\infty) = G/\bar N_{\mathcal{R}},$$ where
$\bar N_{\mathcal{R}}$ is the image of the normal closure $N_{\mathcal{R}}$ of $\mathcal{R}$ in $F$
under the homomorphism $\Phi$.

We will say that $\mathcal{R}_i$ satisfies the condition $(*)$,  if no element $R\in \mathcal{R}_i$ is conjugated  to $R^{-1}$ in $G_{i-1}$.

\section{Asphericity and atoricity}\label{S:aspher}

We will consider diagrams over $G(\infty) = \langle \mathcal{A}\|\mathcal{O}\cup \bigcup_{i=1}^{\infty}\mathcal{R}_i\rangle$  with a gradation similar to \cite{olsh91} (see the definition of a diagram and related notions in \cite{olsh93}). The faces of a diagram are divided into $0$-faces (with labels from $\mathcal{O}$) and $\mathcal{R}$-faces (with labels from $\mathcal{R}=\bigcup_{i=1}^{\infty}\mathcal{R}_i$). The $\mathcal{R}$-faces are subdivided into $\mathcal{R}_i$-faces (with labels from $\mathcal{R}_i$). The rank of a $\mathcal{O}$-face is equal to $0$, the rank of a $\mathcal{R}_i$-face is equal to $i$. Since any diagram contains a finite number of faces, the maximum of the ranks of its faces is well-defined, and will be called {\it the rank of the diagram}. If there are no faces in the diagram, then we assume that its rank is zero.

Following the definition from \cite{olsh93}, a pair of  two different $\mathcal{R}_i$-faces $\Pi_1$ and $\Pi_2$ in a diagram over a  $G$-graded presentation (\ref{infPresent}) will be called  {\it opposite $\mathcal{R}_i$-faces}, if for their boundary labels $R_1$ and $R_2$ reading in a clockwise direction, starting from vertices $o_1$ and $o_2$, one can find (after a series of elementary transformations) a simple path $s$ in the diagram such that $s_{-}=o_1, s_{+}=o_2$ and
$$\varphi(s)^{-1}R_1\varphi(s)R_2=1\,\,\,\text{ in}\,\,\,G_{i-1}.$$

If all $\mathcal{R}_i$ satisfy the condition $(*)$, and for any words $R,R'\in \mathcal{R}_i$ that are not cyclic permutations of each other, $R'$ is not conjugated   to $R$ in $G_{i-1}$, then the word $R_1$ is a cyclic permutation of $R_2^{-1}$.

 A $G$-graded presentation (\ref{infPresent}) will be called {\it $G$-aspherical} ({\it $G$-atorical}), if every diagram on a sphere (on a torus) of positive rank over the presentation (\ref{infPresent}) contains opposite $\mathcal{R}_i$-faces for some $i>0$.

Note that \cite{olsh} uses {\it $i$-pairs of faces} that satisfy stricter requirements than opposite $\mathcal{R}_i$-faces (in the notation above):

1) $R_1$ and $R_2$ are mutually inverse words in the free group;

2) $\varphi(s)=1$ in $G_{i-1}$.

  Similar to \cite{olsh}, a $G$-graded presentation (\ref{infPresent}) will be called {\it aspherical} ({\it atorical}), if any diagram on a sphere (on a torus) of positive rank over the presentation (\ref{infPresent}) contains an $i$-pair of faces for some $i>0$.
Note that if a $G$-graded presentation is aspherical, then all $\mathcal{R}_i$ satisfy the condition $(*)$, and for any words $R,R'\in \mathcal{R}_i$ that are not cyclically permutations of each other, $R'$ is not conjugated  to $R$ in $G_{i-1}$.

It is obvious that asphericity implies $G$-asphericity, and atoricity implies $G$-atoricity.

   As in \cite{olsh91}, let's call  $\tau(\Delta)=(\tau_1,\tau_2,\ldots)$ {\it the type} of a diagram $\Delta$, where $\tau_i$ is the number of faces of rank $i$ in $\Delta$, that is, unlike \cite{olsh}, the type does not take into account the number of $0$-faces in $\Delta$. The types are lexicographically ordered as in \cite{olsh}.

\section{Relations from $[N_{\mathcal{R}}, F]N$}\label{S:commutant2}

In this section, let's consider a group $G(\infty)$ with a presentation (\ref{infPresent}) assuming that all $\mathcal{R}_i$ satisfy the condition $(*)$.

The group $H\cong G/[\bar N_{\mathcal{R}}, G]$ is the central extension of $G(\infty)$. Let's investigate the structure of the central subgroup $\bar
N_R/[\bar N_{\mathcal{R}}, G]$ of the group $H$.

The symmetrized set $\mathcal{R}_i\subset\mathcal{R}$ is divided into disjoint subsets of words conjugated to each other in the group $G_{i-1}$:
 $$\mathcal{R}_i = \bigsqcup_{j} (\mathcal{R}_{i,j}^+\sqcup \mathcal{R}_{i,j}^-),$$ where $\mathcal{R}_{i,j}^- = (\mathcal{R}_{i,j}^+)^{-1}$. Note that $\mathcal{R}_{i,j}^-\neq \mathcal{R}_{i,j}^+$ by the condition $(*)$.
  In each of the sets $\mathcal{R}_{i,j}^+$, fix a representative
${r_{i,j}^+}$.
Let $\mathcal{R}_i^{+}$ denote the set of all representatives $r_{i,j}^+$, and $\mathcal{R}^{+}=\bigcup_{i=1}^{\infty}\mathcal{R}_i^+$.
The group $G_k$, $k=1,2,\ldots$, has the presentation $\langle \mathcal{A}\|\mathcal{O}\cup \bigcup_{i=1}^{k}\mathcal{R}_i^+\rangle$.
The presentation
 \begin{equation} \label{infPresent+}
G(\infty) = \langle \mathcal{A}\|\mathcal{O}\cup \bigcup_{i=1}^{\infty}\mathcal{R}_i^+\rangle
\end{equation}
of the group $G(\infty)$  will be called {\it reduced}.

Consider an arbitrary diagram $\Delta$ over a reduced presentation (\ref{infPresent+}) (which is the same as a diagram over the presentation $\langle \mathcal{A}\|\mathcal{O}\cup \bigcup_{i=1}^{\infty}\mathcal{R}_i^*\rangle$, where $\mathcal{R}_i^*$ is the set containing all elements from $\mathcal{R}_i^{+}$, their inverses and all their cyclic permutations). For each $r_{i,j}^+\in \mathcal{R}_i^{+}$, let $\sigma_{+}(r_{i,j}^+)$ denote the number of the faces in $\Delta$ with labels from $\mathcal{R}_{i,j}^+$, and let $\sigma_{-}(r_{i,j}^+)$ denote the number of the faces with labels from $\mathcal{R}_{i,j}^-$, and $\sigma_{\Delta}(r_{i,j}^+) = \sigma_{+}(r_{i,j}^+) -
\sigma_{-}(r_{i,j}^+)$.
 In the diagram $\Delta$, the number $\sigma_{+}(r_{i,j}^+)$ coincides with the number of faces whose labels are equal to cyclic permutations of the word $r_{i,j}^+$, and $\sigma_{-}(r_{i,j}^+)$ coincides with the number of faces whose labels are equal to cyclic permutations of $(r_{i,j}^+)^{-1}$.

Further, the course of reasoning in Lemma \ref{lemma31.1}, Lemma \ref{lemma31.2} and Theorem \ref{th31.1} coincides with Lemmas 31.1, 31.2 and Theorem 31.1 from \cite{olsh}
(see also \cite{ok3}).
The main difference from \cite{olsh} is the presence of $\mathcal{O}$-faces in diagrams and multipliers from $N$ in words.

\begin{lem}\label{lemma31.1}
If in a disk diagram $\Delta$ over a $G$ graded reduced presentation (\ref{infPresent+}) of a group $G(\infty)$, $\sigma_{\Delta}(R) = 0$
for any $R\in \mathcal{R}^+$, then the boundary label of the diagram $\Delta$ belongs to $[N_{\mathcal{R}}, F]N$.
\end{lem}
\noindent{\it Proof.} In a standard way, making cuts in the  diagram $\Delta$, we get a bouquet of diagrams, each of which contains at most one face. Since the cuts do not change the boundary label as an element of the free group $F$, the boundary label of the diagram $\Delta$ is equal  to
$$(S_1R^{\pm 1
}_1S_1^{-1})\ldots (S_mR^{\pm 1 }_mS_m^{-1})$$
in $F$, where
$R_k \in \mathcal{O}\cup\mathcal{R}^+$, $S_k\in F$.
Each of
these factors belongs to $N_{\mathcal{R}}$ or $N$, therefore, they can
be rearranged modulo $[N_{\mathcal{R}},F]N$. In addition, an element $SRS^{-1}$ is trivial for $R\in \mathcal{O}$ and coincides with $R$ for $R\in \mathcal{R}^+$ modulo $[N_{\mathcal{R}},F]N$. So, by the condition of the lemma, the boundary label is trivial modulo $[N_{\mathcal{R}},F]N$.
\edvo

In order to use induction in the case of $G$-asphericity and $G$-atoricity and prove an analogue of Lemma 31.2 from \cite{olsh}, we need an additional lemma.

\begin{lem}\label{addtolemma31.2} Let a  $G$-graded reduced presentation (\ref{infPresent+}) of a group $G(\infty) = G/\bar N_{\mathcal{R}}$ is $G$-atorical.
Then
\begin{itemize}
\item[$1$)] for an arbitrary diagram $\Delta$ of rank $i$ over
 (\ref{infPresent+}) on a torus, we have $\sigma_{\Delta}(R) = 0$ for any $R\in\mathcal{R}^+$;
\item[$2$)] if $X^{-1}Y^{-1}XY=1$ in $G_i$ for some words $X, Y$,
then there exists a disk diagram $\Delta$ of rank $\leq i$ over the presentation (\ref{infPresent+}) with the boundary label equal to $X^{-1}Y^{-1}XY$,
such that $\sigma_{\Delta}(R) = 0$ for any $R\in\mathcal{R}^+$.
\end{itemize}
\end{lem}

\noindent{\it Proof.} Statements 1) and 2) are proved by common induction on the rank $i$.

1) Let $\Delta$ be a diagram of rank $i$ on a torus. The proof will be by induction on the type $\tau(\Delta)$. The statement is obvious if there are no
$\mathcal{R}$-faces in $\Delta$. Otherwise, as follows from the definition of $G$-atoricity, in $\Delta$, one can find
 a pair of  two different $\mathcal{R}_j$-faces $\Pi_1$ and $\Pi_2$ ($j\leq i$) having boundary labels $R_1$ and $R_2$ reading in a clockwise direction, starting from vertices $o_1$ and $o_2$ and (after a series of elementary transformations) a simple path $s$ such that $s_{-}=o_1, s_{+}=o_2$ and
$\varphi(s)^{-1}R_1\varphi(s)R_2=1\,\,\,\text{ in}\,\,\,G_{j-1},$
moreover, due to the condition $(*)$ and the reducibility of presentation, it can be assumed that $R_1$ and $R_2$ are mutually inverse words in the free group.
For the disk subdiagram $\Gamma$ with boundary label $\varphi(s)^{-1}R_1\varphi(s)R_2$ containing only $\Pi_1$ and $\Pi_2$, it is obvious that $\sigma_{\Gamma}(R) = 0$ for any $R\in\mathcal{R}^+$. By induction assumption,  one can apply statement 2) to the word $\varphi(s)^{-1}R_1\varphi(s)R_2$ and find a disk diagram $\Gamma'$ of rank $\leq j-1$ over the presentation (\ref{infPresent+}) with boundary label $\varphi(s)^{-1}R_1\varphi(s)R_2$, for which $\sigma_{\Gamma'}(R) = 0$ for any $R\in\mathcal{R}^+$.
After removing $\Gamma$ from $\Delta$ and pasting $\Gamma'$, we get a diagram $\Delta'$, for which $\tau(\Delta')<\tau(\Delta)$. By induction assumption, $\sigma_{\Delta'}(R) = 0$ for any $R\in\mathcal{R}^+$. Hence, $\sigma_{\Delta}(R) = \sigma_{\Delta'}(R) - \sigma_{\Gamma'}(R) + \sigma_{\Gamma}(R) = 0$ for any $R\in\mathcal{R}^+$.

2) If $X^{-1}Y^{-1}XY=1$ in $G_i,$
then by the Van Kampen lemma, there exists a disk diagram $\Delta$ of rank $j\leq i$ over the presentation (\ref{infPresent+}) with boundary label equal to $X^{-1}Y^{-1}XY$.
If $j=0$, then it is obvious that $\sigma_{\Delta}(R) = 0$ for any $R\in\mathcal{R}^+$.
If $j>0$, then using the disk diagram $\Delta$, we construct the diagram $\bar{\Delta}$ on a torus in a standard way, gluing the boundary subpath labeled by $X$ with the boundary subpath labeled by $X^{-1}$, and the boundary subpath labeled by $Y$ with the boundary subpath labeled by $Y^{-1}$.
By statement 1)
$\sigma_{\bar\Delta}(R) = 0$ for any $R\in\mathcal{R}^+$.
Hence, $\sigma_{\Delta}(R) = 0$ for any $R\in\mathcal{R}^+$.
 \edvo

\begin{lem}\label{lemma31.2} Let a $G$-graded reduced presentation (\ref{infPresent+}) of a group $G(\infty) = G/\bar N_R$ be $G$-aspherical and $G$-atorical.
Then
\begin{itemize}
\item[$1$)] for an arbitrary spherical diagram $\Delta$ over the presentation
 (\ref{infPresent+}), we have $\sigma_{\Delta}(R) = 0$ for any $R\in\mathcal{R}^+$;
\item[$2$)] if a word $X$ presents in $F$ an element from $[N_{\mathcal{R}}, F]N$,
then for any disk diagram $\Delta$ with boundary label equal to $X$,
$\sigma_{\Delta}(R) = 0$ for any $R\in\mathcal{R}^+$.
\end{itemize}
\end{lem}
\noindent{\it Proof.}
1) Let $\Delta$ be a diagram of rank $i$ on a sphere. The proof will be by induction on the type $\tau(\Delta)$. The statement is obvious if there are no
$\mathcal{R}$-faces in $\Delta$. Otherwise, as follows from the definition of $G$-asphericity,
in $\Delta$, one can find
 a pair of  two different $\mathcal{R}_j$-faces $\Pi_1$ and $\Pi_2$ ($j\leq i$) having boundary labels $R_1$ and $R_2$ reading in a clockwise direction, starting from vertices $o_1$ and $o_2$ and (after a series of elementary transformations) a simple path $s$ such that $s_{-}=o_1, s_{+}=o_2$ and
$\varphi(s)^{-1}R_1\varphi(s)R_2=1\,\,\,\text{ in}\,\,\,G_{j-1},$
moreover, due to the condition $(*)$ and the reducibility of presentation, it can be assumed that $R_1$ and $R_2$ are mutually inverse words in the free group.
For the disk subdiagram $\Gamma$ with boundary label $\varphi(s)^{-1}R_1\varphi(s)R_2$ containing only $\Pi_1$ and $\Pi_2$, it is obvious that $\sigma_{\Gamma}(R) = 0$ for any $R\in\mathcal{R}^+$. By statement 2) of Lemma \ref{addtolemma31.2} there is a disk diagram $\Gamma'$ of rank $\leq j-1$ over the presentation (\ref{infPresent+}) with boundary label $\varphi(s)^{-1}R_1\varphi(s)R_2$, for which $\sigma_{\Gamma'}(R) = 0$ for any $R\in\mathcal{R}^+$.
After removing $\Gamma$ from $\Delta$ and pasting $\Gamma'$, we get a spherical diagram $\Delta'$, for which $\tau(\Delta')<\tau(\Delta)$. By induction assumption, $\sigma_{\Delta'}(R) = 0$ for any $R\in\mathcal{R}^+$. Hence, $\sigma_{\Delta}(R) = \sigma_{\Delta'}(R) - \sigma_{\Gamma'}(R) + \sigma_{\Gamma}(R) = 0$ for any $R\in\mathcal{R}^+$.

2) The word $X$ can be written in $F$ in the form
$(\prod[Y_k,{S_k}R^{\pm 1}_{k}S_k^{-1}])Z,$
where $S_k,Y_k\in F, R_{k}\in \mathcal{R}^+, Z\in N$, since
the relation $[w,uv]=[w,u][uwu^{-1},uvu^{-1}]$ is identical in $F$. With a geometric interpretation of the equality $X=1$ in $G(\infty)$, we obtain a disk diagram $\Delta_0$ with $\sigma_{\Delta_0}(R)=0$ for any $R\in\mathcal{R}^+$. If $\Delta$ is another disk diagram with boundary label equal to
$X$, then $\Delta$ together with the mirror copy $\Delta^0$ of the diagram $\Delta_0$ gives a spherical diagram $\bar \Delta$. By statement 1),
$\sigma_{\bar \Delta}(R) = 0$ for any $R\in\mathcal{R}^+$.
Hence,  $\sigma_{\Delta}(R) = \sigma_{\bar \Delta}(R) - \sigma_{\Delta^0}(R) = 0$ for any $R\in\mathcal{R}^+$. \edvo

\begin{thm}\label{th31.1} Let a $G$-graded reduced presentation (\ref{infPresent+}) of a group $G(\infty) = G/\bar N_\mathcal{R}$ be  $G$-aspherical and $G$-atorical.
Then
\begin{itemize}
\item[$1$)] the following conditions for a word $X\in N_{\mathcal{R}}N$ are equivalent:

a) $X\in [N_{\mathcal{R}}, F]N$;

b) in the writing $Z\prod_k {S_k}R^{\pm 1}_k{S_k^{-1}}$ representing
the word $X$ in $F$, where $Z\in N, R_k\in \mathcal{R}^+, S_k\in F$, the sum of
the exponents of $R$ is equal to zero for each $R\in\mathcal{R}^+$;
\item[$2$)] $\bar N_{\mathcal{R}}/[\bar N_{\mathcal{R}}, G]\simeq N_{\mathcal{R}}N/[N_{\mathcal{R}}, F]N$
is a free Abelian group with basis $\{\bar R\}_{R\in \mathcal{R}^+}$, where $\bar{R}=R[N_{\mathcal{R}}, F]N$.
\end{itemize}
\end{thm}
\noindent{\it Proof.} 1) The condition b) is equivalent to the fact that for a disk diagram $\Delta$ constructed from the writing $Z\prod_k {S_k}R^{\pm 1}_k{S_k^{-1}}$,  for any $R\in\mathcal{R}^+$, $\sigma_{\Delta}(R) = 0$. Therefore, the condition a) follows from b) by statement 2) of Lemma \ref{lemma31.2}, and the condition b) follows from a) by Lemma \ref{lemma31.1}.

2) The group $N_{\mathcal{R}}N$ is generated by all $R\in \mathcal{R}^+$ together with the words conjugated to them in $F$ and their inverses, and by all $Z\in N$. But the words from $N_{\mathcal{R}}$, conjugated in $F$, have the same images in $N_{\mathcal{R}}N/[N_{\mathcal{R}}, F]N$, and the images of the words $Z\in N$ are equal to the identity in $N_{\mathcal{R}}N/[N_{\mathcal{R}}, F]N$. Hence, the Abelian group $N_{\mathcal{R}}N/[N_{\mathcal{R}}, F]N$
is generated by the elements $\{\bar R\}_{R\in \mathcal{R}^+}$.

Let's prove that $\{\bar R\}_{R\in \mathcal{R}^+}$ freely generates $N_{\mathcal{R}}N/[N_{\mathcal{R}}, F]N$. It is enough to make sure that for different elements $R_1,...,R_k\in \mathcal{R}^+$, an equality
$\bar{R}_1^{m_1}...\bar{R}_k^{m_k} = 1$ in $N_{\mathcal{R}}N/[N_{\mathcal{R}}, F]N$ implies the equalities $m_1 = ... = m_k =0$. Indeed, we have the equality
$R_1^{m_1}...R_k^{m_k} = X\in [N_{\mathcal{R}}, F]N$ in the free group. Hence,
$W \equiv R_1^{m_1}...R_k^{m_k}X^{-1}=1$ in $F$.
Applying statement 1) of the theorem to $X$ and to $W$, we get $m_1 =
... = m_k =0$.
 \edvo

 Note that modulo $N$, the full analogue of Lemma 31.2 from \cite{olsh}, in which only asphericity of a $G$-graded reduced presentation is required, can be proved word-by-word as in \cite{olsh}. So in Theorem \ref{th31.1}, it is sufficient for a $G$-graded reduced presentation to be only aspherical.


\section{Use of the limit groups from \cite{olsh91} and \cite{semenov}.}\label{S:commutant1}

Using the limit group constructed in \cite{olsh91}, we obtain analogues of Theorems 31.2, 31.3 
from \cite{olsh} for a non-cyclic torsion-free hyperbolic groups.

As in \cite{olsh91}, let $G^n$ denote the subgroup generated by all $n$-th powers of the elements of a group $G$.

\begin{thm}\label{thm31.2_ex1} For every non-cyclic torsion-free hyperbolic group $G$, there exists an integer $n_0(G)$ such that for every odd $n>n_0(G)$,
the group $H=G/[G^n, G]$ has no torsion, and the subgroup $H^n$ is a free Abelian group.
\end{thm}
\noindent{\it Proof.} Consider the group $G(\infty)=G/G^n$ constructed in \cite{olsh91}.
Then $H^n\cong G^n/[G^n, G]$, and $H/H^n$ is the group $G/G^n$. As proved in \cite{olsh91} (p. 565), for $G(\infty)$, the analogue of Lemma 18.2 from \cite{olsh} holds, that is, there are asphericity and atoricity, so one can apply Theorem \ref{th31.1}, i.e. the analogue of Theorem 31.1 from \cite{olsh}. In addition, as proved in \cite{olsh91} (p. 565),
for $G(\infty)$, the analogue of Lemma 18.3 from \cite{olsh} holds, that is, any word is conjugated to a power of a period of some rank in $H/H^n$.
  Taking into account the presence of these analogues, the proof of Theorem \ref{thm31.2_ex1} repeats the proof of Theorem 31.2 from \cite{olsh}.
 \edvo

\begin{thm}\label{thm31.3_ex1}
For every non-cyclic torsion-free hyperbolic group, there exists a non-Abelian torsion-free quotient group in which any two non-trivial subgroups have a non-trivial intersection.
\end{thm}
\noindent{\it Proof.} For the group $H$ from Theorem \ref{thm31.2_ex1}, the reasoning from the proof of Theorem 31.3 \cite{olsh} is repeated, only the reference to Theorem 31.2 in \cite{olsh} should be replaced with the reference to its analogue -- to Theorem \ref{thm31.2_ex1}.
 \edvo


Using the limit group constructed in \cite{semenov, semenov_phd}, we obtain generalizations of Theorems 31.4, 31.7 and 31.8 from \cite{olsh}.

 \begin{thm}\label{thm31.4_ex2} For every non-cyclic torsion-free hyperbolic group, there exists a non-Abelian torsion-free quotient group in which
 any two non-trivial subgroups have a non-trivial intersection and all proper subgroups are cyclic.
\end{thm}
\noindent{\it Proof.} Consider an arbitrary non-cyclic torsion-free hyperbolic group $G_0$. Since the group $G_0$ is not elementary, its commutant $G_0'$ is not elementary.  Since $G_0$ is a non-cyclic torsion-free hyperbolic group, then by Proposition 1, Theorems 1, 2 of \cite{olsh93}, there exists a  homomorphism from the group $G_0$ to a non-cyclic torsion-free hyperbolic group $G$ such that it is surjective on $G_0'$.

 By Theorem 2.6 \cite{semenov_phd} (= Theorem 1 \cite{semenov}), for  the group $G$ and for any odd $n$ greater than some sufficiently large number $n_0=n_0(G)$, there exists an infinite quotient group
 $G(\infty)=G/\bar N_{\mathcal{R}}$, all proper subgroups of which are cyclic groups of order dividing $n$.
By Lemma 54 and Theorem 2.4 in \cite{semenov_phd}, the analogues of Lemma 25.1 and Theorem 26.4 from \cite{olsh} are true for $G(\infty)$. Therefore, we can continue as in the proof of Theorem 31.3 \cite{olsh}.

Consider the group $H=G/[\bar N_{\mathcal{R}}, G]\simeq F/[N_{\mathcal{R}}, F]N$ (hereinafter the designations from Sections \ref{S:df} and \ref{S:commutant2}). Since there are asphericity and atoricity (Lemma 54 \cite{semenov_phd}), Theorem \ref{th31.1} can be applied. Hence, the Abelian group
$\bar N_{\mathcal{R}}/[\bar N_{\mathcal{R}}, G]\simeq N_{\mathcal{R}}N/[N_{\mathcal{R}}, F]N$ freely generated by the elements $\{\bar{R}\}_{R\in \mathcal{R}^+}$. In $H$, consider the subgroup $L$ consisted of the products $\prod_k \bar{R}_k^{s_k}$, where $\sum_ks_k=0$. Since $\bar{R}$ lies in the center of $H$, the subgroup $L$ is normal in $H$. In the quotient group $\mathbb{A}=H/L$, the coset $\bar{R}_1L=\bar{R}_2L=\ldots=C$ has infinite order. The group $\mathbb{A}$ is an extension of the group $G(\infty)$ by the infinite cyclic central subgroup $\langle C\rangle$. Since $G(\infty)$ is not Abelian and infinite, then $\mathbb{A}$ is also not Abelian and infinite.

To prove that there are no elements of finite order in $\mathbb{A}$, assume the opposite. For every nontrivial element $X$ of finite order in $\mathbb{A}$, we have
$X\notin \langle C\rangle$, i.e. $X\neq 1$ in $G(\infty)$, and there exists $j\geq 1$ such that $X$ has finite order in $G_j$.  By Theorem 2.4 \cite{semenov_phd}, $X$ is conjugated in $G_j$ to a power of a period $A$ of the first type, i.e. $A^n\in\mathcal{R}$. Hence, replacing $X$ with a conjugate, we can assume that $X=A^kC^l$ in the group $\mathbb{A}$.
 Therefore, $X^n=A^{kn}C^{ln}=C^kC^{ln}$. Since the element $C$ has infinite order in $\mathbb{A}$, we get that $n$ divides $k$, that contradicts the condition $X\notin \langle C\rangle$.

Since there is no nontrivial element of finite order in $\mathbb{A}$, and all elements in $\mathbb{A}/\langle C\rangle$ have finite orders, any  element of the group $\mathbb{A}$ in some power lies in $\langle C\rangle$. So the intersection of any two nontrivial subgroups of the group $\mathbb{A}$ is nontrivial.

Consider an arbitrary proper subgroup $K$ of the group $\mathbb{A}$. It is mapped to a subgroup $\hat{K}$ of the group $G(\infty)$. The subgroup $\hat{K}$
is either finite, or coincides with $G(\infty)$. In the first case, $K$ is elementary, and therefore it is cyclic. In the second case, $\mathbb{A}=K\langle C\rangle$, hence, $K$ is a normal subgroup of the group $\mathbb{A}$.
Since $\mathbb{A}/K$ is an Abelian group, the commutant of the group $\mathbb{A}$ belongs to $K$.
Since the commutant of the group $G$ coincides with $G$, the commutant of $\mathbb{A}$ coincides with $\mathbb{A}$. Hence, $K=\mathbb{A}$, i.e. all proper subgroups of the group $\mathbb{A}$ are cyclic.
 \edvo

 If the orders of all elements of a group are bounded in common,  the smallest common multiple of the orders of all its elements will be called the {\it period} of the group.
  A periodic group will be called a {\it $p$-group}, if the order of each of its elements is some power of a prime number $p$. A $p$-subgroup of an arbitrary group will be called a {\it Sylow $p$-subgroup} if it is not contained in a larger $p$-subgroup of this group.

 \begin{thm}\label{thm31.7_ex1}  For every prime $p$, for every non-cyclic torsion-free hyperbolic group, there exists a non-Abelian quotient group $K$ of finite period, in which all proper subgroups are cyclic, and the Sylow $p$-subgroup $P$ is central, but it is not a direct multiplier of the group $K$.
\end{thm}
\noindent{\it Proof.} Consider the group $\mathbb{A}$ from the Theorem \ref{thm31.4_ex2} for a sufficiently large odd number $n$, not divisible by $p$. Since the infinite cyclic subgroup $\langle C\rangle$ is central, the subgroup $\langle C^p\rangle$ is normal in $\mathbb{A}$. Put $K=\mathbb{A}/\langle C^p\rangle$. Since the group $G(\infty)$ is simple (Remark 2 to Theorem 2.6 of \cite{semenov_phd}), the center of the group $K$ coincides with the subgroup $P=\langle C\rangle/\langle C^p\rangle$. Since the orders of the elements of the group $G(\infty)=K/P$ are finite and divide $n$, then $P$ is a Sylow $p$-subgroup in $K$, and the orders of the elements of the group $K$ are finite and divide $np$. Since the complete preimage of any subgroup of the group $K$ is a subgroup in  $\mathbb{A}$, then all proper subgroups of the group $K$ are cyclic. In particular, it follows that $P$ is not a direct multiplier of the group $K$.
  \edvo

  Note that in the Theorem \ref{thm31.7_ex1}, any proper subgroup of the group $K$ either trivially intersects with $P$, or coincides with $P$, or contains $P$ as a direct multiplier. And additionally (compared to Theorem 31.7 of \cite{olsh}), there is a cyclicity of all proper subgroups in $K$.

 \begin{thm}\label{thm31.8_ex1} For every torsion-free non-cyclic hyperbolic group, for every sufficiently large prime number $p$ and any integer $k\geq 0$, there exists a non-Abelian quotient group $K$ of the period $p^{k+1}$, which contains a central cyclic subgroup $D$ of order $p^k$ such that any subgroup of the group $K$  is cyclic (except $K$) and either contains $D$, or is contained in $D$.
 \end{thm}
\noindent{\it Proof} is similar to the proof of Theorem 31.8 \cite{olsh}. As a group $S$ from the proof of Theorem 31.8 \cite{olsh}, we need to take the group  $\mathbb{A}$  from Theorem \ref{thm31.4_ex2} for $n=p$. Instead of referring to Theorem 26.4 from \cite{olsh}, we need to refer to its analogue -- to Theorem 2.4 from \cite{semenov_phd}.
 \edvo

\section{Use of a limit group from \cite{olsh93}.}\label{S:ex2}

In \cite{olsh93} (and subsequently in another way in \cite{semenov}), for any non-cyclic torsion-free hyperbolic group, A.Yu. Ol'shanskii constructed a non-Abelian torsion-free quotient group $\overline{G}$, all of whose proper subgroups are cyclic, and the intersection of any two different maximal proper subgroups is trivial. (The latter follows from Lemma 10 \cite{olsh91} and the fact that $\overline{G}$ is a limit of torsion-free hyperbolic groups.)
To obtain some examples of central extensions of $\overline{G}$,  we need the following theorem. For completeness, we formulate and prove it not only for the group $\overline{G}$ from Corollary 1 \cite{olsh93}, but also for $\hat{G}$ from Corollary 4 \cite{olsh93} and  $\tilde{G}$ from Corollary 2 \cite{olsh93}.


\begin{thm}\label{thm31.2} Let $G$ be an arbitrary non-cyclic torsion-free hyperbolic group, $G/\bar{N}_{\mathcal{R}}$ be the group $\hat{G}$ constructed in Corollary 4 \cite{olsh93}, or $\overline{G}$ from Corollary 1 \cite{olsh93}, or $\tilde{G}$ from Corollary 2 \cite{olsh93} (we assume that
  $\overline{G}$ and $\tilde{G}$ are constructed using Theorem 4 \cite{olsh93} instead of Theorem 2 \cite{olsh93} or Theorem 3 \cite{olsh93}). Then $H=G/[\bar N_{\mathcal{R}}, G]$ has no torsion, and
$\bar N_{\mathcal{R}}/[\bar N_{\mathcal{R}}, G]$ is a free Abelian group with countable basis of elements $\{\bar R\}_{R\in \mathcal{R}^+}$, where $\bar{R}=R[N_{\mathcal{R}}, F]N$.
 \end{thm}
\noindent{\it Proof.} Consider the group $\hat{G}$
(for $\overline{G}$ and $\tilde{G}$, everything is similar, even simpler).

Recall that in Corollary 4 \cite{olsh93},   $\hat{G}$ is an infinite non-Abelian group  such that all its proper subgroups are finite. In constructing of $\hat{G}= G(\infty)= {\rm lim}_{\rightarrow}\,\, G_i$, one consider all elements $g_1,g_2,\ldots$ of $G$ and all its finitely generated non-elementary subgroups $H_1, H_2,\ldots$, and by Theorem 4 \cite{olsh93}, one find a sequence of epimorphisms $G_0=G\rightarrow G_1\rightarrow G_2\rightarrow\ldots$ such that every element $g_i$ has  finite order in $G_{2i-1}$ (as in Corollary 2 \cite{olsh93}) and the image of $H_i$ is elementary in $G_{2i}$ or coincides with $G_{2i}$ (as in Corollary 1 \cite{olsh93}).
For sufficiently large numbers $m_{i,0}$,  the numbers $m_i$ are chosen arbitrarily so that $m_i\geq m_{i,0}$, and the numbers $m_{2i-1}$ are odd.
The sets $\mathcal{R}_{2i}$ are the sets $\mathcal{R}_{k_{2i},l_{2i},m_{2i}}$ from the proof of Theorem 2 \cite{olsh93} or empty sets, if $H_i$ is elementary in $G_{2i-1}$.
If $g_i$ has finite order in $G_{2i-2}$, then $G_{2i-2}\rightarrow G_{2i-1}$ is an identical homomorphism and $\mathcal{R}_{2i-1}=\emptyset$.
Otherwise, one consider the maximal elementary subgroup $E_{i,0}$ containing $g_i$, and its maximal
infinite cyclic normal subgroup $\langle \bar{g}_i\rangle$, and $G_{2i-1}=G_{2i-2}/\langle \overline{g}_i^{m_{2i-1}} \rangle^{G_{2i-2}}$, and the set  $\mathcal{R}_{2i-1}$ is given as the set of all cyclic permutations of the word $W_i^{\pm m_{2i-1}}$, where $W_i$ is a word representing the element $\bar{g}_i$.

Since $G$ is a non-cyclic torsion-free hyperbolic group, then, as shown in the proof of Corollary 4 \cite{olsh93}, $G$ is non-elementary and satisfies
the quasiidentical relation
$$(\forall x,y\in G) x^2y=yx^2\Rightarrow xy=yx.$$ 
In addition, by construction (Theorem 4 \cite{olsh93}), all $G_i$ are non-elementary hyperbolic groups satisfying the quasiidentical relation
(for this, $m_{2i-1}$ are chosen odd).

Since it follows from $xyx^{-1}=y^{-1}$ that $x^2yx^{-2}=y$, then by the quasiidentical relation, we obtain $y^2=1$, which is impossible for an element of infinite order.
Since all elements of $\mathcal{R}_i$ have infinite order in $G_{i-1}$ (by Lemmas 4.1(1), 4.2, 7.2 \cite{olsh93}), for any $R\in \mathcal{R}_i$, $R^{-1}$ is not conjugated to $R$ in $G_{i-1}$. Thus, $\mathcal{R}_i$ satisfy the $C_2(\epsilon_i,\mu_i,\lambda_i,c_i,\rho_i)$-condition by Lemma 4.1 (2) \cite{olsh93} for odd $i$ and the $C_1(\epsilon_i,\mu_i,\lambda_i,c_i,\rho_i)$-condition by Lemma 4.2 \cite{olsh93} for even $i$ (there exist appropriate $\lambda_i>0, (l_i>0), c_i\geq 0$ such that for any $\mu_i>0, \epsilon_i\geq 0, \rho_i>0$, there exists $m_{i,0}>0$).

Since the $C_1(\epsilon_i,\mu_i,\lambda_i,c_i,\rho_i)$-condition and the $C_2(\epsilon_i,\mu_i,\lambda_i,c_i,\rho_i)$-condition hold, by Lemma 6.6 of \cite{olsh93}, the presentation (\ref{infPresent}) of the group $\hat{G}=G(\infty)$ is $G$-aspherical, and also, as noted in the proof of Theorem 4 in \cite{olsh93}, it is $G$-atorical.

Consider the group $H=G/[\bar N_{\mathcal{R}}, G]\simeq F/[N_{\mathcal{R}}, F]N$ (the designations from Sections \ref{S:df} and \ref{S:commutant2}). By Theorem  \ref{th31.1}, the Abelian group
$M=\bar N_{\mathcal{R}}/[\bar N_{\mathcal{R}}, G]\simeq N_{\mathcal{R}}N/[N_{\mathcal{R}}, F]N$ is freely generated by the elements of $\{\bar{R}\}_{R\in \mathcal{R}^+}$,  where $\bar{R}=R[N_{\mathcal{R}}, F]N$.
Since $\hat{G}$ is not hyperbolic, $\mathcal{R}^+\neq\cup_{i=1}^{k}\mathcal{R}_i^+$ for every $k$ and the set $\mathcal{R}^+=\cup_{i=1}^{\infty}\mathcal{R}_i^+$ is infinite. The group $H$ is an extension of the group $\hat{G}$ by the central subgroup $M$.

Let's prove that there is no torsion in $H$. Assume the opposite. For every nontrivial element $X$ of finite order in $H$, we have that
$X\notin M$, i.e. $X\neq 1$ in $\hat{G}$, and there exists the smallest number $j$ such that $X$ has finite order in $G_j$, and $j\geq 1$.
Among all such elements, we choose an element $X\in H$ with the smallest number $j$. By Lemma 7.2 \cite{olsh93}, in $G_j$, $X$ is conjugated to some element represented by a word  $A$, which represents an element from the centralizer $C_{G_{i-1}}(R_i)$ for some $R_i\in \mathcal{R}_i^+$, $i\leq j$.
Therefore, since $G_{i-1}$ is a hyperbolic group, there exists a natural number $s$ such that $A^s=R_i^l$ in $G_{i-1}$,
so,  $A^s=R_i^l\prod S_kR_{i_k}^{\delta_k}S_k^{-1}$ in $G$, where $\delta_k=\pm 1, R_{i_k}\in \bigcup_{\nu=1}^{i-1}\mathcal{R}_{\nu}^+$.
 In $H$, we have
 $$A^s=\bar{R}_i^l\prod_{\nu=1}^{i-1}\prod_{k_\nu}\bar{R}_{k_{\nu}}^{l_{k_{\nu}}},$$
  where $R_{k_{\nu}}\in \mathcal{R}_{\nu}^+$ are different for different ${k_{\nu}}$, $l_{k_{\nu}}$ are integers. Replacing $X$ with a conjugate, we can assume that $X=AZ$ in the group $H$, where $Z\in M$.
 Since $Z$ belongs to the center, in the group $H$, we have
   $$X^s=A^sZ^{s}=\bar{R}_i^l\prod_{\nu=1}^{i-1}\prod_{k_\nu}\bar{R}_{k_{\nu}}^{l_{k_{\nu}}}Z^{s}.$$
   Since $\bar{R}_i$, $\bar{R}_{k_{\nu}}$ (for all $k_{\nu}$ for all $\nu=1,\ldots,i-1$) are different basis elements of $M$, and there are no elements of finite order in $M$, then $X^s=1$ in $H$ and $l=sd, l_{k_{\nu}}=sd_{k_{\nu}}$ for all $k_{\nu}$.

 Hence, $(A(\bar{R}_i^{d}\prod_{\nu=1}^{i-1}\prod_{k_{\nu}}\bar{R}_{k_{\nu}}^{d_{k_{\nu}}})^{-1})^s=1$ in $H$, and $A(\bar{R}_i^{d}\prod_{\nu=1}^{i-1}\prod_{k_{\nu}}\bar{R}_{k_{\nu}}^{d_{k_{\nu}}})^{-1}\neq 1$ in $H$ (otherwise it would be $A\in M$, which means $X\in M$). Since $A\in C_{G_{i-1}}(R_i)$, then  we have $(A(\bar{R}_i^{d}\prod_{\nu=1}^{i-1}\prod_{k_{\nu}}\bar{R}_{k_{\nu}}^{d_{k_{\nu}}})^{-1})^s=1$ in $G_{i-1}$. Since $i-1<j$, we have a contradiction with the choice of $X$.
 \edvo

 The proof of the absence of torsion in $H$ for $H$ of Theorem \ref{thm31.2} is a little different from one for $H$ of Theorem 31.2\cite{olsh}, since $l$ does not have to equal to $1$ here.

\begin{thm}\label{thm31.4_ex3} For every non-cyclic torsion-free hyperbolic group and every $k\in \mathbb{N}$, there exists a non-Abelian torsion-free quotient group $\mathbf{A}$, in which all maximal proper subgroups are isomorphic to $\mathbb{Z}^k$ and the intersection of any two different maximal proper subgroups coincides with the center of $\mathbf{A}$, which is isomorphic to $\mathbb{Z}^{k-1}$, and no nontrivial subgroup of the center is a direct multiplier of $\mathbf{A}$.
\end{thm}
\noindent{\it Proof.} Consider an arbitrary non-cyclic torsion-free hyperbolic group. As in the proof of Theorem \ref{thm31.4_ex2}, it is surjectively mapped to a non-cyclic torsion-free hyperbolic group $G$, whose commutant coincides with $G$.

  By Corollary 1  of \cite{olsh93}, for the group $G$, there exists a non-Abelian torsion-free quotient group $G(\infty)=\overline{G}=G/\bar N_{\mathcal{R}}$ such that all its proper subgroups are cyclic, the intersection of any two different maximal subgroups is trivial, and the center is trivial. This group proves the statement for $k=1$. Let $k>1$.

Consider the group $H=G/[\bar N_{\mathcal{R}}, G]$. By Theorem \ref{thm31.2}, the Abelian group
$\bar N_{\mathcal{R}}/[\bar N_{\mathcal{R}}, G]$  has a countable basis of free generators $\{\bar{R}\}_{R\in \mathcal{R}^+}$. Numerate the elements of $\{\bar{R}\}_{R\in \mathcal{R}^+}$ by the natural numbers.
Since $\bar{R}$ lies in the center of $H$, the subgroup $L$ consisted with the products $\prod_m \bar{R}_{i_m}^{s_m}$, where $\sum_ms_m=0$, $i_m\geq k-1$, is normal in $H$. The group $\mathbf{A}=H/L$ is an extension of the group $\overline{G}$ by the free Abelian central subgroup $\mathbf{C}=\langle \bar{R}_1L\rangle\times\ldots\times\langle \bar{R}_{k-1}L\rangle$. Since $\overline{G}$ is non-Abelian and infinite, then $\mathbf{A}$ is also non-Abelian and infinite. Since the group $\overline{G}$ has no torsion, then the group $\mathbf{A}$ also has no torsion.
Since the center of $\overline{G}$ is trivial, the center $Z(\mathbf{A})$ of $\mathbf{A}$ coincides with $\mathbf{C}$.

Consider an arbitrary proper subgroup $K$ of $\mathbf{A}$. It is mapped to a subgroup $\hat{K}$ of $\overline{G}$. Either the subgroup $\hat{K}$ coincides with $\overline{G}$, or it is a cyclic subgroup generated by some element $\hat{g}$. In the first case, $\mathbf{A}=KZ(\mathbf{A})$, and hence, $K$ is a normal subgroup of $\mathbf{A}$. Since $\mathbf{A}/K$ is an Abelian group, the commutant of $\mathbf{A}$ belongs to $K$, and hence, $K=\mathbf{A}$ (it follows similarly that no nontrivial subgroup of the center is a direct multiplier of $\mathbf{A}$).  In the second case, if $\hat{g}$ is not trivial, then  $K=\langle g\rangle\times (K\cap Z(\mathbf{A}))\leq\langle g\rangle\times Z(\mathbf{A})\neq \mathbf{A}$, where $g\in K$ is a preimage of $\hat{g}$, otherwise $K\leq Z(\mathbf{A})\leq\langle g\rangle\times Z(\mathbf{A})\neq \mathbf{A}$ for an arbitrary element $g\notin Z(\mathbf{A})$.
 \edvo

 \begin{thm}\label{thm31.4_ex4} For every non-cyclic torsion-free hyperbolic group, there exists a non-Abelian torsion-free quotient group $\mathbf{A}$ such that all its maximal proper subgroups are free Abelian of countable rank and the intersection of any two different maximal proper subgroups coincides with the center of $\mathbf{A}$, which is a free Abelian group with basis $\{e_i|i\in \mathbb{N}\}$, and a basis of any maximal proper subgroup of $\mathbf{A}$ can be obtained by adding one element to $\{e_i|i\in \mathbb{N}\}$, and no nontrivial subgroup of the center is a direct multiplier of $\mathbf{A}$.
\end{thm}

\noindent{\it Proof} is similar to the proof of Theorem \ref{thm31.4_ex3} for trivial $L$ and $\mathbf{C}=\bar N_{\mathcal{R}}/[\bar N_{\mathcal{R}}, G]$. \edvo

Note that the center of $\mathbf{A}$ from Theorems  \ref{thm31.4_ex3} (from Theorem \ref{thm31.4_ex4}) is the Frattini subgroup.

\begin{thm}\label{thm31.4_ex5}
For every non-cyclic torsion-free hyperbolic group, there exists an infinite non-Abelian quotient group $K$ such that all its elements of finite order form a subgroup $P$, which is the center, and  $P$ is not a direct multiplier of $K$.
\end{thm}
\noindent{\it Proof.} Consider the group $\mathbf{A}$ from Theorem \ref{thm31.4_ex3}, for example, for $k=2$.
In this case, its center $\mathbf{C}$ is an infinite cyclic subgroup $\langle C\rangle$. For an arbitrary integer $s>1$, its subgroup $\langle C^s\rangle$ is normal in $\mathbf{A}$. Put $K=\mathbf{A}/\langle C^s\rangle$.
Since the center of $\overline{G}=K/P$ is trivial, the center of $K$ coincides with the subgroup $P=\langle C\rangle/\langle C^s\rangle$.
Since $\overline{G}$ has no torsion, all elements of finite order in $K$ lie in $P$. The subgroup $P$ is not a direct multiplier of $K$, since the complete preimage of any proper subgroup of $K$ is a proper subgroup of $\mathbf{A}$, and $K$ is non-Abelian.
 \edvo

{\it Keywords:} Torsion-free hyperbolic groups, quotients


\end{document}